\documentclass[a4paper,10pt]{amsart}
\usepackage{amsmath,amssymb,amsfonts,amsthm}
\usepackage{latexsym,mathrsfs}
\usepackage{graphicx}
\usepackage[all]{xy}
\input xypic
\usepackage{color}
\usepackage{hyperref}
\hypersetup{colorlinks=true,linkcolor=red,citecolor=blue}
\usepackage[T1]{fontenc}
\usepackage[utf8]{inputenc}
\usepackage{lmodern}
\usepackage{microtype}
\usepackage{enumerate}
\theoremstyle{plain}
\newtheorem{theorem}[subsection]{{\bf Theorem}}
\newtheorem*{theorem*}{{\bf Theorem}}

\newtheorem*{corollary*}{{\bf Corollary}}
\newtheorem{proposition}[subsection]{{\bf Proposition}}
\newtheorem{lemma}[subsection]{{\bf Lemma}}

\newtheorem{conjecture}[subsection]{{\bf Conjecture}}
\theoremstyle{definition}

\theoremstyle{remark}

\numberwithin{equation}{subsection}

\DeclareMathOperator{\HH}{H}
\DeclareMathOperator{\K}{K}
\DeclareMathOperator{\B}{B}

\DeclareMathOperator{\M}{M}

\newcommand{\QZ}{\mathbb{Q}/\mathbb{Z}}

\newcommand{\cwedge}{\curlywedge}

\begin{document}
\title[Exponents of Bogomolov multipliers]{On the exponent of Bogomolov multipliers}
\author{Primo\v z Moravec}
\address{{
Faculty of  Mathematics and Physics, University of Ljubljana,
and Institute of Mathematics, Physics and Mechanics,
Slovenia}}
\email{primoz.moravec@fmf.uni-lj.si}
\subjclass[2010]{14E08, 20F12, 20F45}
\keywords{Bogomolov multiplier, exponent, finite group}
\thanks{The author acknowledges the financial support from the Slovenian Research Agency (research core funding No. P1-0222, and projects No. J1-8132, J1-7256 and N1-0061).}
\date{\today}
\begin{abstract}
\noindent
We prove that if $G$ is a finite group, then the exponent of its Bogomolov multiplier divides the exponent of $G$ in the following four cases: (i) $G$ is metabelian, (ii) $\exp G=4$, (iii) $G$ is nilpotent of class $\le 5$, or (iv) $G$ is a $4$-Engel group.
\end{abstract}
\maketitle
\section{Introduction}
\label{s:intro}
\noindent
Let $G$ be a group given by a free presentation $G=F/R$. Then the {\it Schur multiplier} of $G$ can be defined via Hopf's formula as
$\M (G)=([F,F]\cap R)/[F,R]$. It is known that if $G$ is a finite group, then $\M (G)\cong \HH_2(G,\mathbb{Z})\cong \HH^2(G,\mathbb{Q}/\mathbb{Z})$; see, for example, Beyl and Tappe's book \cite{Bey82} for details and applications of Schur multipliers.

A significant part of the theory of Schur multipliers consists of estimating their size, rank, or exponent. Focusing on the latter one, it was already Schur \cite{Sch07} who showed that if $G$ is a finite group, then $(\exp \M (G))^2$ divides $|G|$. On the other hand, it was shown in \cite{Mor07} that, for every finite group $G$, the exponent of $\M(G)$ can be bounded in terms of $\exp G$ only. Good bounds of this kind are however still out of reach. In practice it often happens that $\exp \M(G)$ divides $\exp G$. One may conjecture that this is always the case, yet a counterexample of exponent 4 with Schur multiplier of exponent 8 was constructed by Bayes, Kautsky and Wamsley \cite{Bay73}. We mention that there no known examples of odd order groups $G$ with $\exp \M (G)>\exp G$, though it seems plausible that they exist.
On the other hand, $\exp \M (G)$ always divides $\exp G$ at least when $G$ is of one of the following types:
nilpotent of class $\le 3$ \cite{Jon74,Mor08a}, powerful $p$-group \cite{Lub87}, a $p$-group of maximal class \cite{Mor11}. In addition to that, we have the following.
\begin{theorem}[\cite{Mor07,Mor08a,Mor08}]
\label{t:expm}
Let $G$ be a group of finite exponent.
\begin{enumerate}
  \item If $G$ is nilpotent of class $\le 4$, then $\exp \M(G)$ divides $2\cdot \exp G$.
  \item If $G$ is a $4$-Engel group, then $\exp \M(G)$ divides $10\cdot \exp G$.
  \item If $\exp G=4$, then $\exp \M(G)$ divides 8.
  \item If $G$ is metabelian, then $\exp \M(G)$ divides $(\exp G)^2$.
\end{enumerate}
\end{theorem}
There are numerous other estimates for $\exp \M(G)$, but we do not go further into this here, see, for example, Sambonet's paper \cite{Sam15} for a short review of such results.

Given a free presentation $G=F/R$, define
\begin{align*}
\B_0(G) &= ([F,F]\cap R) / (\langle \K(F)\cap R\rangle),\\
G\cwedge G &= [F,F] / (\langle \K(F)\cap R\rangle),
\end{align*}
where $\K(F)$ is the set of commutators in $F$.
It is shown in \cite{Mor12} that if $G$ is a finite group and $V$ a faithful representation of $G$ over $\mathbb{C}$, then the dual of $\B_0(G)$ is naturally isomorphic to the unramified Brauer group $\HH^2_{\rm nr}(\mathbb{C}(V)^G,\QZ )$
introduced by Artin and Mumford \cite{Art72}. This invariant represents an obstruction to {\it Noether's problem} \cite{Noe16} asking as to whether the field $\mathbb{C}(V)^G$ is purely transcendental over $\mathbb{C}$. The above mentioned result of \cite{Mor12} is based on a result of Bogomolov \cite{Bog88} who showed that $\HH^2_{\rm nr}(\mathbb{C}(V)^G,\QZ )$ is naturally isomorphic to the intersection of the kernels of restriction maps $\HH^2(G,\QZ)\to \HH^2(A,\QZ)$,
where $A$ runs through all (two-generator) abelian subgroups of $G$. The latter group is also known as the {\it Bogomolov multiplier} of $G$. Here we use the same name for $\B_0(G)$.

Our main result is the following.

\begin{theorem}
  \label{t:expb0}
  Let $G$ be a group of finite exponent. If $G$ satisfies one of the following properties, then $\exp (G\curlywedge G)$ divides $\exp G$:
  \begin{enumerate}
    \item Nilpotent of class $\le 5$,
    \item Metabelian,
    \item 4-Engel,
    \item $\exp G=4$.
  \end{enumerate}
\end{theorem}

This result thus complements Theorem \ref{t:expm}. It also implies, that, in the cases listed in Theorem \ref{t:expb0}, it always happens that $\exp \B_0(G)$ divides $\exp G$. Note that the question of determining $\exp \B_0(G)$ was recently addresed by Garc\' ia-Rodr\' iguez, Jaikin-Zapirain, and Jezernik
\cite[Theorem 6]{Gar17}.

In view of the above result and further extensive computational evidence we pose the following conjecture:

\begin{conjecture}
  \label{conj:expb0}
  Let $G$ be a finite group. Then $\exp\B_0(G)$ divides $\exp G$, and thus, in particular, $\exp \M(G)$ divides $(\exp G)^2$.
\end{conjecture}

In order to justify this conjecture, we first mention that Saltman \cite{Sal84} showed that for every $p>2$ there exists a $p$-group of exponent $p$ with non-trivial Bogomolov multiplier of exponent $p$, thus the bound given in the conjecture is sharp. On the other hand, if $\exp G$ is not prime, then $\B_0(G)$ is usually small (or of small exponent), compared to $G$. Fernandez-Alcober and Jezernik \cite{Fer18} recently constructed examples of finite $p$-groups $G$ of maximal class with $\exp\B_0(G)\approx\sqrt{\exp G}$, and it appears that this might be close to the worst case. In fact, we have not been able to find a group $G$ with $\exp (G\curlywedge G)>\exp G$, so it may be possible that even the
stronger conjecture with $\B_0(G)$ replaced by $G\curlywedge G$ might hold true.

The proof of Theorem \ref{t:expb0} goes roughly as follows. First, one may assume without loss of generality that $G$ is a finite group. By \cite{Jez15a}, there exists a so-called CP cover $H$ of $G$ (see Section \ref{s:prelim} for the details) whose main feature is that $G\curlywedge G$ is isomorphic to $[H,H]$, and all commutator relations of $G$ lift to commutator relations in $H$. The calculations are then performed in $[H,H]$. In the metabelian case, the proof is fairly straightforward. In the remaining cases, it relies on a careful examination of the power structure of the lower central series of $H$, which uses information on the free groups in a given variety. For class $\le 5$ groups, this can be obtained using elementary commutator calculus. In the exponent 4 and 4-Engel case, we mainly use M. Hall's description of the free 3-generator group of exponent 4 \cite{Hal73}, and Nickel's computations of free 4-Engel groups of low ranks \cite{Nic99}, along with Havas and Vaughan-Lee's proof of local nilpotency of 4-Engel groups \cite{Hav05}.

\section{Preliminaries}
\label{s:prelim}

\subsection{CP covers}
\label{ss:cp}

Let $G$ be a group and $Z$ a $G$-module. Denote by $e=(\chi , H,\pi )$ the extension
\begin{equation*}
\xymatrix{ 1\ar[r] & Z\ar[r]^\chi & H\ar[r]^\pi & G\ar[r] & 1}
\end{equation*}
of $Z$ by $G$.
Following \cite{Mor12}, we say that $e$ is a {\it CP extension} if commuting
pairs of elements of $G$ have commuting lifts in $H$.
A stem central CP extension of $Z$ by $G$, where $|Z|=|\B_0(G)|$, is called a {\it CP cover} of $G$. CP covers are analogs of the usual covers in the theory of Schur multipliers. It is proved in \cite{Jez15a} that every finite group has a CP cover. It also follows from \cite{Jez15a} that if $e=(\chi , H,\pi )$ of the above form is a CP cover of $G$, then $Z^\chi\cap \K (H)=1$. This in particular implies that any commutator law satisfied by $G$ is also satisfied by $H$.

\subsection{Collection process}
\label{ss:collect}

\noindent
We recall \cite[Theorem 11.2.4]{Hal59} that Hall's Basis Theorem implies that if $F$ is a free nilpotent group of class $c$ and $a,b\in F$, then the word $(ab)^n$, where $n$ is a non-negative integer, can be written uniquely as a product $c_1^{n_1}c_2^{n_2}\cdots c_t^{n_t}$, where $c_i$ are basic commutators in $\{a,b\}$ of weights $1,2,\ldots ,c$, and $n_i=b_1{n\choose 1}+b_2{n\choose 2}+\cdots +b_r{n\choose r}$,
where $r$ is the weight of $c_i$ and $b_j$ are non-negative integers not depending on $n$. Specifically, we will be interested in the case when $F$ is free nilpotent of class 6. We need to determine the coefficients $b_i$ explicitly, and this can be done using the collection process described in \cite[Section 12.3]{Hal59}. We omit the details regarding calculations, and only record the values of $b_i$ for all basic commutators of weight $\le 6$ in $\{ a,b\}$ in Table \ref{tablebi}.

\begin{table}[h!tb]
\begin{tabular}{|l|c|c|c|c|c|c|}
\hline
Commutator $c_i$ & $b_1$ & $b_2$ & $b_3$ & $b_4$ & $b_5$ & $b_6$\\
\hline
$a$ & 1 & & & & & \\
$b$ & 1 & & & & & \\
$[b,a]$ & & 1 & & & & \\
$[b,a,a]$ & & & 1 & & & \\
$[b,a,b]$ & & 1 & 2 & & & \\
$[b,a,a,a]$ & & & & 1 & & \\
$[b,a,a,b]$ & & & 2 & 3 & & \\
$[b,a,b,b]$ & & & 2 & 3 & & \\
$[b,a,a,[b,a]]$ & & & 1 & 7 & 6 & \\
$[b,a,b,[b,a]]$ & & & 6 & 18 & 12 & \\
$[b,a,a,a,a]$ & & & & & 1 & \\
$[b,a,a,a,b]$ & & & & 3 & 4 & \\
$[b,a,a,b,b]$ & & & 1 & 6 & 6 & \\
$[b,a,b,b,b]$ & & & & 3 & 4 & \\
$[b,a,b,[b,a,a]]$ & & & 4 & 21 & 36 & 20\\
$[b,a,a,a,[b,a]]$ & & & & 3 & 13 & 10\\
$[b,a,a,b,[b,a]]$ & & & 2 & 24 & 52 & 30\\
$[b,a,b,b,[b,a]]$ & & & 3 & 27 & 54 & 30\\
$[b,a,a,a,a,a]$ & & & & & & 1\\
$[b,a,a,a,a,b]$ & & & & & 4 & 5\\
$[b,a,a,a,b,b]$ & & & & 3 & 12 & 10\\
$[b,a,a,b,b,b]$ & & & & 3 & 12 & 10\\
$[b,a,b,b,b,b]$ & & & & & 4 & 5\\
\hline
\end{tabular}
\caption{Coefficients in exponents of $c_i$.}
\label{tablebi}
\end{table}


%


\section{Proof of Theorem \ref{t:expb0}}
\label{s:proof}

\noindent
In what follows, $G$ will be a group of finite exponent satisfying one of the properties listed in Theorem \ref{t:expb0}. In each of those cases, $G$ is then locally finite. As $\B_0$ commutes with direct limits \cite{Mor12}, one may assume without loss of generality that $G$ is a finite group; furthermore, Bogomolov's results \cite{Bog88} imply that we can restrict ourselves to the case when $G$ is a finite $p$-group. Let
$$\xymatrix{ 1\ar[r] & Z\ar[r] & H \ar[r]^{\pi} & G\ar[r] & 1}$$
be a CP cover of $G$, where $Z$ is a central subgroup of $H$ with the property that $Z\cong \B_0(G)$ and $Z\cap \K(H)=1$. From here on the proof goes on by considering each case separately.

\subsection{Metabelian groups}
\label{s:meta}

The case of metabelian groups is easy:

\begin{theorem}
\label{t:metabel}
Let $G$ be a metabelian group of finite exponent. Then the exponent of $G\cwedge G$ divides $\exp G$.
\end{theorem}

\proof
Put $\exp G=e$. Note that $H$ is also metabelian, hence it suffices to prove that $[x,y]^e=1$
for all $x,y\in H$. We expand
$1=[x,y^e]=[x,y]^e\prod _{k=2}^e[x,{}_ky]^{{e\choose k}}$. Observe that $\prod _{k=2}^e[x,{}_ky]^{{e\choose k}}\in Z$. Furthermore,
$$\prod _{k=2}^e[x,{}_ky]^{{e\choose k}}=\left [ \prod _{k=2}^e[x,{}_{k-1}y]^{{e\choose k}},y\right ]\in\K(H),$$ therefore $[x,y]^e=1$, as required.
\endproof


\subsection{Exponent 4}
\label{s:exp4}

At first we list some properties of groups of exponent 4 that will be used in the proof of this case.

\begin{lemma}
\label{l:exp4}
Let $G$ be a group of exponent 4 and $a,b,c\in G$.
\begin{enumerate}[(a)]
\item The group $\langle a,b\rangle$ is nilpotent of class $\le 5$, $\langle a,b,c\rangle$ is nilpotent of class $\le 7$, and $\langle [a,b],c\rangle$ is nilpotent of class $\le 4$,
\item $[[a,b]^2,a]=1$,
\item $[a,b,a,a^2[a,b]]=1$,
\item $[c,[a,b],[a,b],[a,b]]=1$,
\end{enumerate}
\end{lemma}

\proof
All the above properties can be deduced immediately from a polycyclic presentation of $B(3,4)$, see also \cite{Hal73}.
\endproof

\begin{theorem}
\label{t:exp4}
Let $G$ be a group of exponent 4. Then the exponent of $G\cwedge G$ divides 4.
\end{theorem}

\proof
As noted above, we may assume without loss of generality that $G$ is a finite group.
Choose $x,y,z\in H$.

First note that $[[x,y]^2,x]\in Z\cap \K(H)=1$ by Lemma \ref{l:exp4}, therefore
\begin{equation}
\label{eq:exp4_1}
1=[x,y,x]^2[x,y,x,[x,y]].
\end{equation}

Take $w\in \{ x,y\}$. As $\langle x,y\rangle$ is nilpotent of class $\le 5$, we get
$1=[[x,y]^2,x,w]=[x,y,x,w]^2$. From here it follows that
\begin{equation}
\label{eq:exp4_2}
\gamma _4(\langle x,y\rangle )^2=1.
\end{equation}
We also have that
$[x,y,z]^4=1$ by \cite[Proof of Theorem 2.6]{Mor07}.
Now we expand $1=[x^4,y]$ using \cite[Lemma 9]{Mor08a}:

\begin{align*}
1 &= [x^4,y]\\
  &= [x,y]^4[x,y,x]^6[x,y,x,x]^4[x,y,x,x,x][x,y,x,[x,y]]^{14}\\
  &= [x,y]^4[x,y,x]^2[x,y,x,x,x].
\end{align*}

Lemma \ref{l:exp4} implies that $[x,y,x,x^2[x,y]]=1$. Expanding this using the class restriction, we obtain $1=[x,y,x,x]^2[x,y,x,x,x][x,y,x,[x,y]]$, and this implies
$[x,y,x,x,x]=[x,y,x,[x,y]]$. From \eqref{eq:exp4_1} and the above expansion we get that $[x,y]^4=1$.

By Lemma \ref{l:exp4}, the group $\langle [x,y],z\rangle$ is metabelian and nilpotent of class $\le 4$, and we also have that $[z,[x,y],[x,y],[x,y]]=1$. We now expand $([x,y]z)^4$ using Subsection \ref{ss:collect} and \eqref{eq:exp4_2}:
$$([x,y]z)^4=z^4[z,[x,y]]^2[z,[x,y],[x,y],z][z,[x,y],z,z].$$
Denote $w=[z,[x,y]]^2[z,[x,y],[x,y],z][z,[x,y],z,z]$ and consider the following
words:
\begin{align*}
  w_1 &= [x,z,z,x^2z^2[z,y,x][z,x,y,y]],\\
  w_2 &= [y^2z^2,[z,y,z]],\\
  w_3 &= [y^2z^2,[y,x,x,x][z,x,x,x]],\\
  w_4 &= [y^2z^2,[z,x,z,z][z,y,y,x,x]],\\
  w_5 &= [[z,y][z,x]z^2,z^2[z,x][z,y][z,x,x][z,x,y][z,y,x][z,y,y][z,x,x,x][z,y,y,x,x]],\\
  w_6 &= [[z,y][z,x]z^2, [z,y,z,x][z,y,z,y][z,y,z,x,x]],\\
  w_7 &= [[z,y][z,x]z^2, [z,x,z,z]].
\end{align*}
The subgroup $\langle x,y,z\rangle$ is an image of
$$K=\langle a,b,c \mid \hbox{class } 7, \hbox{laws } [x_1^4,x_2]=[x_1,x_2]^4=[[x_1,x_2]^2,x_1]=[x_1,_{3}[x_2,x_3]]=1\rangle.$$
Expansions of the above defined words in $K$ into products of basic commutators reveals that $w=w_1w_2\cdots w_7$. On the other hand, inspection of the presentation of $B(3,4)$ shows that $w_i\in \K(H)\cap Z=1$ for all $i=1,2,\ldots ,7$, therefore $w=1$.
This immediately implies $([x,y]z)^4=z^4$ for all $x,y,z\in  H$. From here it
is not difficult to conclude that $\exp \gamma _2(H)$ divides 4, and this finishes the proof.
\endproof


\subsection{4-Engel groups}
\label{s:4engel}

\noindent
The aim of this section is to prove
\begin{theorem}
\label{t:4eng}
Let $G$ be a $4$-Engel group of finite exponent. Then the exponent of $G\cwedge G$ divides $\exp G$.
\end{theorem}

As $4$-Engel groups are locally nilpotent \cite{Hav05}, the situation can be easily reduced to the case when $G$ is a finite $p$-group. If $p\neq 2,5$,
then it follows from \cite{Mor08} that even $\exp (G\wedge G)$ divides $\exp G$. This is no longer true when $p=2$ or $p=5$. In the case when $p=5$, there is a short proof of
Theorem \ref{t:4eng}. Let $H$ be a CP cover of $G$ and denote $\exp G=5^e$. Note that $H$ is a 4-Engel 5-group, hence it is regular \cite{Hav07}. It follows from \cite{Mor08} that if $x,y\in H$, then $[x,y]^{5^e}=1$. Regularity now implies that $\gamma_2(H)$ has exponent dividing $5^e$.

We are thus left with $4$-Engel $2$-groups. The argument here is more involved. We start with some preliminaries.

\begin{lemma}
\label{l:4engexp}
Let $G$ be a $4$-Engel group of exponent $2^e$ and $a,b,c\in G$.
\begin{enumerate}[(a)]
\item $\gamma _7(\langle a,b\rangle )=\gamma _8(\langle a,b,c\rangle )^2=\gamma _9(\langle a,b,c\rangle )=1$,
\item $[a,b,a]^{2^{e-1}}=[a,b,b]^{2^{e-1}}=1$,
\item $\gamma _4(\langle a,b,c\rangle )^{2^{e-1}}=1$,
\end{enumerate}
\end{lemma}

\proof
It follows from \cite{Nic99} that if $\langle a,b,c\rangle$ is a $4$-Engel group, then $\gamma _7(\langle a,b\rangle )=(\gamma _8(\langle a,b,c\rangle )/\gamma _9(\langle a,b,c\rangle ))^{30}=\gamma _9(\langle a,b,c\rangle )^3=1$. This proves (a). The fact that the exponent of $\gamma _4(\langle a,b,c\rangle )$ divides $2^{e-1}$ is proved in \cite[Lemma 4.6]{Mor08}, whereas
the proof of that Lemma also yields (b).
\endproof

We will also use the following:

\begin{lemma}[cf Lemma 4.4 of \cite{Mor08}]
\label{l:xny}
Let $G$ be a 4-Engel group, $a,b\in G$ and $n$ a non-negative integer. Then
\begin{equation*}
[a^n,b]=
[a,b]^n[a,b,a]^{{n\choose 2}}[a,b,a,a]^{{n\choose 3}}[a,b,a,[a,b]]^{{n\choose 2}+2{n\choose 3}}
\end{equation*}
\end{lemma}

Referring to a polycyclic presentation of the free 4-Engel group with two or three generators obtained in \cite{Nic99}, we have:

\begin{lemma}
\label{l:e24}
Let $G$ be a $4$-Engel group and $a,b\in G$.
\begin{enumerate}[(a)]
\item $[b,a,a,[b,a],a]=[b,a,a,a,b,a]=1$,
\item $[b,a,a,b,a]^3[b,a,b,b,a]=[a,b,a,[a,b]][b,a,b,a,b,a]^3$,
\item If $G$ has no elements of order 3, then $[c,[a,b],[a,b],[a,b]]\in\gamma _7(\langle a,b,c\rangle )^2\gamma _8(\langle a,b,c\rangle )$.
\end{enumerate}
\end{lemma}

\begin{proposition}
\label{p:4engid}
Let $G$ be 4-Engel group of exponent $2^e$. Then
\begin{enumerate}[(a)]
\item $[[a,b]^{2^{e-1}},a]=1$,
\item $[c,[a,b],[a,b],[a,b]]^{2^{e-2}}=1$.
\end{enumerate}
\end{proposition}

\proof
We may assume that $e>2$.
Let us expand $(ab)^{2^e}=1$ using Subsection \ref{ss:collect} and Lemma \ref{l:4engexp}. We obtain
\begin{equation}
\label{eq:4engid_0}
\begin{split}
[a,b]^{{2^e \choose 2}} ={} & ( [b,a,a,a][b,a,a,b][b,a,b,b][b,a,a,[b,a]][b,a,a,a,b]\\
& \times [b,a,b,[b,a,a]][b,a,a,a,b,b] ) ^{{2^e\choose 4}}.
\end{split}
\end{equation}
We commute this with $a$ and apply class restriction and Lemma \ref{l:e24} (a):
\begin{equation}
\label{eq:4engid_1}
[[a,b]^{{2^e\choose 2}},a]=\left ( [b,a,a,b,a][b,a,b,b,a]\right ) ^{{2^e \choose 4}}.
\end{equation}
Using Lemma \ref{l:4engexp} and Lemma \ref{l:xny}, we obtain after a short calculation that
\begin{equation}
\label{eq:4engid_2}
[[a,b]^{2^{e-1}},a]=[a,b,a,[a,b]]^{{2^{e-1}\choose 2}}.
\end{equation}
The equations \eqref{eq:4engid_1} and \eqref{eq:4engid_2}, together with Lemma \ref{l:e24} (b), give
$
[b,a,b,a,b,a]^{2^{e-2}}=1.
$
This immediately yields $\gamma _6(\langle a,b\rangle )^{2^{e-2}}=1$. Now replace $b$ by $ab$ in
\eqref{eq:4engid_0} and use \eqref{eq:4engid_0}. Expansion under given class restriction gives
\begin{equation}
\label{eq:4engid_3}
1=\left ( [b,a,b,a][b,a,a,b]\right )^{2^{e-2}}.
\end{equation}
If we further replace $b$ by $ab$ in \eqref{eq:4engid_3} and apply \eqref{eq:4engid_3}, we obtain
$[b,a,a,b,a]^{2^{e-2}}=1$. Replacing $a$ by $ba$ in this identity, we conclude that also
$[b,a,b,b,a]^{2^{e-2}}=1$. Equation \eqref{eq:4engid_1} now gives $[[a,b]^{{2^e\choose 2}},a]=1$. This proves (a), whereas (b) follows directly from
Lemma \ref{l:e24} (c) and Lemma \ref{l:4engexp} (c), as $e>2$.
\endproof

\begin{theorem}
\label{t:4eng2e}
Let $G$ be a $4$-Engel group of exponent $2^e$. Then the exponent of $G\cwedge G$ divides $2^e$.
\end{theorem}

\proof
Note that $H$ is a 4-Engel group. Take $x,y,z\in H$ and let $a=x^\pi$, $b=y^\pi$, $c=z^\pi$. Proposition \ref{p:4engid} implies
\begin{align*}
1 &= [[x,y]^{2^{e-1}},x]\\
 &= [x,y,x]^{2^{e-1}}[x,y,x,[x,y]]^{2^{e-1}\choose 2}.
\end{align*}
Equation \eqref{eq:4engid_2} implies $[a,b,a,[a,b]]^{2^{e-1}\choose 2}=1$, therefore
$[x,y,x,[x,y]]^{2^{e-1}\choose 2}=[[x,y,x]^{2^{e-1}\choose 2},[x,y]]\in \K (H)\cap Z=1$. This gives $[x,y,x]^{2^{e-1}}=1$. From Lemma \ref{l:xny} we get
$$1=[x^{2^e},y]=[x,y]^{2^e}[x,y,x]^{2^e\choose 2}[x,y,x,x]^{2^e \choose 3}[x,y,x,[x,y]]^{{2^e\choose 2}+2{2^e\choose 3}}.$$
 Using the above equations, we see that this identity
implies $[x,y]^{2^e}=1$. Now note that the subgroup $\langle [x,y],z\rangle$ is nilpotent of class $\le 5$, since $H$ is 4-Engel. We expand $([x,y]z)^{2^e}$ using collection process (see Subsection \ref{ss:collect}):
\begin{equation}
\label{eq:4engH_1}
([x,y]z)^{2^e}=z^{2^e}([z,[x,y],[x,y],[x,y]][z,[x,y],[x,y],z][z,[x,y],z,z])^{2^e\choose 4}.
\end{equation}
Note that $[z,[x,y],[x,y],[x,y]]^{2^e\choose 4}=[[z,[x,y]]^{2^e\choose 4},[x,y],[x,y]]\in \K(H)$, and Proposition \ref{p:4engid} implies that $[z,[x,y],[x,y],[x,y]]^{2^e\choose 4}\in Z$. This immediately shows that
$[z,[x,y],[x,y],[x,y]]^{2^e\choose 4}=1$. Thus $([z,[x,y],[x,y],z][z,[x,y],z,z])^{2^e\choose 4}\in Z$. Furthermore, the class restriction yields that
$([z,[x,y],[x,y],z][z,[x,y],z,z])^{2^e\choose 4}=[[z,[x,y],[x,y]][z,[x,y],z],z]^{2^e\choose 4}=[([z,[x,y],[x,y]][z,[x,y],z])^{2^e\choose 4},z]\in \K(H)$, therefore we conclude that
$([z,[x,y],[x,y],z][z,[x,y],z,z])^{2^e\choose 4}=1$. Equation \eqref{eq:4engH_1} thus gives $([x,y]z)^{2^e}=z^{2^e}$, and induction on the commutator length shows that $\exp H'$ divides $2^e$.
\endproof


\subsection{Groups of nilpotency class $\le 5$}
\label{s:nilpotent}

We will prove the following:

\begin{theorem}
\label{t:class5}
Let $G$ be a group of finite exponent and class $\le 5$. Then the exponent of $G\cwedge G$ divides $\exp G$.
\end{theorem}

Again, we may assume that $G$ is a finite $p$-group of class $\le 5$ and exponent $p^e$. The CP cover $H$ of $G$ is then also nilpotent of class $\le 5$. Let $x,y\in H$. Assume first that $p>2$. Then it follows from the proof of \cite[Theorem 13]{Mor08a} that $[x,y]^{p^e}=1$. As $[H,H]$ is nilpotent of class $\le 2$, it is regular, hence $\exp [H,H]$ divides $p^e$. Thus we are left with
the case when $G$ is a $2$-group. Without loss we can asssume that $e>2$.

Take $g,h,k\in G$. Then the expansion of $(g[h,k])^{2^e}=1$ yields
$$1=[h,k,g]^{{2^e\choose 2}}[h,k,g,[h,k]]^{{2^e\choose 2}+2{2^e\choose 3}}
[h,k,g,g,g]^{{2^e\choose 4}}.$$
If we replace $h$ by a commutator $[h_1,h_2]$ in the above equation, we get, after renaming the variables, that
$$[h_1,h_2,h_3,h_4]^{2^{e-1}}=1,$$
therefore $\gamma _4(G)^{2^{e-1}}=1$. By the class restriction this implies
$\gamma _4(H)^{2^{e-1}}=1$. Take now $x,y,z\in H$. Then
$$1=[[x,y]^{2^e},z]=[x,y,z]^{2^e}[x,y,z,[x,y]]^{{2^e\choose 2}}=[x,y,z]^{2^e},$$
hence $\gamma _3(H)^{2^e}=1$. As $[H,H]$ is nilpotent of class $\le 2$, it suffices to prove that $[x,y]^{2^e}=1$ for all $x,y\in H$, and then Theorem \ref{t:class5} follows.

Take $x,y\in H$. Then
\begin{equation}
  \label{eq:c5e1}
  1=[x^{2^e},y]=[x,y]^{2^e}[x,y,x]^{{2^e \choose 2}}[x,y,x,x,x]^{{2^e\choose 4}}.
\end{equation}
If we interchange $x$ and $y$ in \eqref{eq:c5e1}, we get
\begin{equation}
  \label{eq:c5e2}
  1=[x,y]^{2^e}[x,y,y]^{{2^e\choose 2}}[x,y,y,y,y]^{{2^e\choose 4}}.
\end{equation}
Now we replace $x$ by $yx$ in \eqref{eq:c5e1} and apply \eqref{eq:c5e1} and \eqref{eq:c5e2}. After a short calculation we obtain
\begin{equation}
  \label{eq:c5e3}
  [x,y]^{2^e}=\left ( [x,y,x,x,y][x,y,x,y,x][x,y,x,y,y][x,y,y,x,x][x,y,y,x,y][x,y,y,y,x]\right )^{{2^e\choose 4}}.
\end{equation}
As $H$ is nilpotent of class $\le 5$, we have that
$[x,y,x,y,x]=[x,y,y,x,x]$ and $[x,y,x,y,y]=[x,y,y,x,y]$. Thus \eqref{eq:c5e3} can be rewritten as
\begin{equation}
  \label{eq:c5e4}
  [x,y]^{2^e}=\left ( [x,y,x,x,y][x,y,y,y,x]\right )^{{2^e\choose 4}}.
\end{equation}
Denote $f={2^e\choose 4}$, and
\begin{align*}
  u &= [y,x][y,x,y],\\
  v &= [y,x]^{-1}[y,x,y]^{-1}[y,x,x,x][y,x,y,x],\\
  w &= [y^fu, y^{-f}v].
\end{align*}
We expand $w$:
\begin{align*}
  w &= [y^fu,v][y^fu,y^{-f}]^v\\
  &= [y^f,v]^u[u,v][u,y^{-f}]^v\\
  &= [y^f,v][y^f,v,u][u,v][u,y^{-f}][u,y^{-f},v]\\
  &= [y^f,v][u,y^{-f}][u,v]\left ([y,v,u][u,y,v]^{-1} \right )^f.
\end{align*}
Note that
  $$[u,v] =[[y,x],[y,x,y]^{-1}][[y,x,y],[y,x]^{-1}]=[[y,x],[y,x,y]]^{-1}[[y,x,y],[y,x]]^{-1}=1,$$
  and the Hall-Witt identity gives
  $1=[y,v,u][v,u,y][u,y,v]=[y,v,u][u,y,v]$.
  Thus $w=[y^f,v][u,y^{-f}]$.
  Now,
  $$[u,y^{-f}]=[(y^{-1})^f,u]^{-1}=[y^{-1},u]^{-f}[y^{-1},u,y^{-1}]^{-{f\choose 2}}[y^{-1},u,y^{-1},y^{-1}]^{-{f\choose 3}}.$$
  Quick calculation shows that $[y^{-1},u]=[y,x,y]$,
  $[y^{-1},u,y^{-1}]=[y,x,y,y,y][y,x,y,y]^{-1}$, and $[y^{-1},u,y^{-1},y^{-1}]=[y,x,y,y,y]$. Therefore
  $$[u,y^{-f}]=[y,x,y]^{-f}[y,x,y,y]^{{f\choose 2}}[y,x,y,y,y]^{-{f\choose 2}-{f\choose 3}}.$$
  On the other hand, we easily get that $$[y,v]=[y,x,y][y,x,y,[x,y]][y,x,y,y][y,x,x,x,y][y,x,y,x,y],$$ and thus
  \begin{align*}
    [y^f,v] &= [y,v]^f[y,v,y]^{{f\choose 2}}[y,v,y,y]^{{f\choose 3}}\\
    &= [y,x,y]^f[y,x,y,[x,y]]^f[y,x,y,y]^{f+{f\choose 2}}[y,x,x,x,y]^f
    [y,x,y,x,y]^f[y,x,y,y,y]^{{f\choose 2}+{f\choose 3}}.
  \end{align*}
  We thus get
  $$w=[y,x,y,y]^{f+2{f\choose 2}}[y,x,y,[x,y]]^f\left ([y,x,x,x,y][y,x,y,x,y]\right ) ^f.$$
  Note that $[y,x,y,y]^{f+2{f\choose 2}}=[y,x,y,y]^{f^2}=1$, since $f^2$ is divisible by $2^{2e-4}\ge 2^{e-1}$. As $H$ has class $\le 5$, we also have
  $[y,x,y,[x,y]]=[y,x,y,x,y][x,y,y,y,x]$. This, together with \eqref{eq:c5e4}, implies
  $w=\left ([x,y,x,x,y][x,y,y,y,x]\right )^f=[x,y]^{2^e}$. We conclude
  that $[x,y]^{2^e}\in \K(H)\cap Z=1$, and this proves Theorem \ref{t:class5}.




\end{document}